\documentclass{article}
\usepackage{latexsym}
\usepackage{amssymb}
\begin{document}

% \mmbox enables macros to survive outside of $ ... $
\newcommand{\mmbox}[1]{\mbox{${#1}$}}
\newcommand{\proj}[1]{\mmbox{{\bf P}^{#1}}}
\newcommand{\affine}[1]{\mmbox{{\bf A}^{#1}}}
\newcommand{\Ann}[1]{\mmbox{{\rm Ann}({#1})}}
\newcommand{\caps}[3]{\mmbox{{#1}_{#2} \cap \ldots \cap {#1}_{#3}}}
\newcommand{\N}{{\bf N}}
\newcommand{\Z}{{\bf Z}}
\newcommand{\R}{{\bf R}}
\newcommand{\Tor}{\mathop{\rm Tor}\nolimits}
\newcommand{\Ext}{\mathop{\rm Ext}\nolimits}
\newcommand{\Hom}{\mathop{\rm Hom}\nolimits}
\newcommand{\im}{\mathop{\rm Im}\nolimits}
\newcommand{\rank}{\mathop{\rm rank}\nolimits}
\newcommand{\supp}{\mathop{\rm supp}\nolimits}
\newcommand{\text}[1]{\mbox{\rm {#1}}}

\sloppy
\newtheorem{defn0}{Definition}[section]
\newtheorem{prop0}[defn0]{Proposition}
\newtheorem{conj0}[defn0]{Conjecture}
\newtheorem{thm0}[defn0]{Theorem}
\newtheorem{lem0}[defn0]{Lemma}
\newtheorem{corollary0}[defn0]{Corollary}
\newtheorem{example0}[defn0]{Example}
\newtheorem{rem0}[defn0]{Remark}

\newenvironment{defn}{\begin{defn0}}{\end{defn0}}
\newenvironment{prop}{\begin{prop0}}{\end{prop0}}
\newenvironment{conj}{\begin{conj0}}{\end{conj0}}
\newenvironment{thm}{\begin{thm0}}{\end{thm0}}
\newenvironment{lem}{\begin{lem0}}{\end{lem0}}
\newenvironment{cor}{\begin{corollary0}}{\end{corollary0}}
\newenvironment{exm}{\begin{example0}\rm}{\end{example0}}
\newenvironment{rem}{\begin{rem0}\rm}{\end{rem0}}

\newcommand{\defref}[1]{Definition~\ref{#1}}
\newcommand{\propref}[1]{Proposition~\ref{#1}}
\newcommand{\thmref}[1]{Theorem~\ref{#1}}
\newcommand{\lemref}[1]{Lemma~\ref{#1}}
\newcommand{\corref}[1]{Corollary~\ref{#1}}
\newcommand{\exref}[1]{Example~\ref{#1}}
\newcommand{\secref}[1]{Section~\ref{#1}}
\newcommand{\remref}[1]{Remark~\ref{#1}}

\newcommand{\qed}{\mbox{$\Box$}}
\newenvironment{proof}{\noindent {\bf Proof.}}{\qed\vskip 6pt}
\newenvironment{proofn}{\noindent {\bf Proof.}}{\vskip 6pt}

\newcommand{\Macaulay}{{\it Macaulay\/}}
\newcommand{\std}{Gr\"{o}bner}
\newcommand{\Std}{Gr\"{o}bner}
\newcommand{\poina}{\pi({\cal A}, t)}
\newcommand{\poinaa}{\overline{\pi}({\cal A}, t)}
\newcommand{\chia}{\chi({\cal A},t)}
\newcommand{\poinap}{\pi({\cal A}', t)}
\newcommand{\poinapp}{\pi({\cal A}'', t)}
\newcommand{\tm}{D_0({\cal A})}
\newcommand{\ts}{{\cal D}_0({\cal A})}
\newcommand{\tmp}{D_0({\cal A}')}
\newcommand{\tsp}{{\cal D}_0({\cal A}')}
\newcommand{\tmpp}{D_0({\cal A}'')}
\newcommand{\tspp}{{\cal D}_0({\cal A}'')}
\newcommand{\Dp}{D^p({\cal A})}
\newcommand{\Do}{D^1({\cal A})}
\newcommand{\TS}{{\cal D}({\cal A})}
\newcommand{\rjq}{R/J_Q}
\newcommand{\rjqm}{R/<J_Q,x>}
\newcommand{\rjp}{R/J_{Q'}}
\newcommand{\rjqx}{R(-1)/<J_Q:x>}
\newcommand{\rjqxx}{R(-2)/<J_Q:x^2>}
\newcommand{\rjqxm}{R(-1)/<J_Q:x,x>}
\newcommand{\jqmp}{<J_Q:x^2>/J_{Q'}}
\newcommand{\jqxx}{<J_Q:x^2>}
\newcommand{\jqx}{<J_Q:x>}
\newcommand{\jqxm}{<J_Q:x,x>}
\newcommand{\jp}{J_{Q'}}
\newcommand{\jq}{J_{Q}}
\newcommand{\ooa}{\Omega^1({\cal A})}
\newcommand{\oo}{\Omega^1_0}
\newcommand{\opa}{\Omega^p({\cal A})}
\newcommand{\wpa}{\Lambda^p (\Omega^1({\cal A}))}
\newcommand{\wps}{\Lambda^p({\rm Der}\,_k S)}
\newcommand{\wpdo}{\Lambda^p(D^1)}
\newcommand{\wpdoo}{\Lambda^p(D^1_0)}
\newcommand{\wpmdoo}{\Lambda^{p-1}(D^1_0)}
\newcommand{\wpdos}{\Lambda^p(\widetilde{D^1})}
\newcommand{\logforms}{\Omega^1(\mbox{log  }D)}
\parskip = 4pt

\title {The module of logarithmic p-forms of a locally free arrangement}
\author {
        Mircea Musta\c t\v a\thanks{Corresponding author} and
        Hal Schenck\thanks{Partially supported by an NSF postdoctoral research
fellowship} \\
        $^*$University of California, Berkeley\\
        $^\dagger$Northeastern University\\
        mustata@math.berkeley.edu\\
        schenck@neu.edu}

\maketitle
\sloppy

\begin{abstract} \noindent For an essential, central hyperplane arrangement
${\cal A}
\subseteq V \simeq k^{n+1}$ we show that $\ooa$ (the module of logarithmic one
forms with poles along ${\cal A}$) gives rise
to a locally free sheaf on ${\bf P}^n$ if and only if for
all $X \in L_{\cal A}$ with rank $X  <$ dim $V$, the module $\Omega^1({\cal
A}_X)$ is free. Our main result says that in this case
 $\poina$ is essentially the
Chern polynomial. The proof is based on a result of Solomon-Terao
\cite{st} and a formula we give for the Chern polynomial of a bundle
${\cal E}$ on ${\bf P}^n$ in terms of the Hilbert series of
$\oplus_{m\in{\bf Z}}H^0({\bf P}^n,\wedge^i{\cal E}(m))$.
If $\ooa$ has projective
dimension one and is locally free, we give a minimal free resolution for
$\Omega^p$, and show that
 $\Lambda^p \ooa \simeq \Omega^p({\cal A})$, generalizing
results of Rose-Terao on generic arrangements.
\end{abstract}

\maketitle

\pagebreak

\section{Introduction}\label{sec:intro}

If $X$ is a complex manifold and $D$ a divisor with normal crossings
($D = \sum D_i$, $D_i$ smooth and meet transversely), then
associated to $D$ is the sheaf $\logforms$ of meromorphic one forms
with logarithmic poles on $D$. Deligne introduced this sheaf
in \cite{d} and shows (among other things) that it is locally free.
Dolgachev and Kapranov seem to have been the first to examine
in depth the case where $D$ is a set of hyperplanes in general position
in ${\bf P}^n$; one striking result they obtain is that if $D = \cup_{i=1}^{d}
H_i$ and $d \ge 2n+3$, then the hyperplanes can be recovered from $\logforms$
unless the $H_i$ osculate a rational normal curve of degree $n$.

We also consider the case when the divisor is a set
of hyperplanes in ${\bf P}^n$, but assume only that the
hyperplanes are distinct. There are two main themes of this
paper. In \cite{st}, Solomon and Terao give a formula for the
Poincar\'e polynomial of an (essential, central) arrangement
in terms of the Hilbert series of certain graded modules $D^i$
associated to the arrangement. The formula generalizes Terao's
famous freeness theorem \cite{t}: If $D^1$ is a free module, then
the Poincar\'e polynomial factors. This suggests a connection to
Chern polynomials; motivated by the Solomon-Terao result we prove
a formula relating the Chern polynomial of a bundle ${\cal E}$ on
${\bf P}^n$ to the Hilbert series of the modules
$\oplus_{m\in{\bf Z}}H^0({\bf P}^n,\wedge^i{\cal E}(m))$.

Since an arbitrary arrangement does not have normal crossings,
$\logforms$ is in general no longer locally
free. Silvotti \cite{si} studies this situation, and remedies the
problem by blowing up the arrangement at the non-normal crossings;
$\sigma^{*}D$
automatically has normal crossings on the blowup X, so yields a locally
free sheaf
$\Omega^1_X(\mbox{log }\sigma^{*}D) \simeq \oplus {\cal F}_j$.
Using a vanishing result of Esnault, Schechtman and  Viehweg \cite{esv},
Silvotti obtains a formula for the  coefficients of the Poincar\'e polynomial
in terms of $\chi(\Lambda^i{\cal F}_j)$. However, the computations can
be quite complicated; in particular, Silvotti does not recover
Terao's theorem.

The second point of this paper is that even when the hyperplanes
are not in general position, there are situations where $\logforms$
is a vector bundle on ${\bf P}^n$. We relate $\logforms$ to the
modules $D^i$ mentioned above, and prove a criterion for the associated
sheaves to be locally free. This class of arrangements was studied by
Yuzvinsky in \cite{y2}; Yuzvinsky proves that for such arrangements
the Hilbert polynomial of $D^1$ is a combinatorial invariant. We show
that the Chern polynomial of the dual of $D^1$
is in fact the Poincar\'e polynomial of the arrangement (truncated by
$t^{n+1}$). Hence, the Hilbert polynomial of $D^1$ may be obtained from
the Poincar\'e polynomial via Hirzebruch-Riemann-Roch. We close with
some results specific to the situation where $\Omega^1$ or $D^1$ has
projective dimension one. First, we review some facts about arrangements.

\section{Hyperplane Arrangements}\label{sec:ha}
A hyperplane arrangement ${\mathcal A}$ is a finite collection of
codimension one
linear subspaces of a fixed vector space V. ${\mathcal A}$ is {\it central}
if each hyperplane contains the origin {\bf 0} of V. A fundamental invariant
of ${\mathcal A}$ is the Poincar\'e polynomial $\poina$. There are
various ways  of defining $\poina$; the simplest is from the intersection
lattice
$L_{\mathcal A}$ of ${\mathcal A}$. $L_{\mathcal A}$ consists of the
intersections of
the elements of ${\mathcal A}$, ordered by reverse inclusion.
The rank function on $L_{\cal A}$ is given by the codimension in $V$. V is the
lattice element $\hat{0}$; the rank one elements are the hyperplanes
themselves. ${\mathcal A}$ is called {\it essential} if
rank $L_{\mathcal A} =$ dim $V$. Henceforth, unless explicitly
stated otherwise,
 all arrangements will be
{\it essential} and {\it central}, and $V$ will be $k^{n+1}$, with $k$
an arbitrary field.
We briefly review some fundamental definitions; for more information see
Orlik and Terao
(\cite{ot}).

\begin{defn}
The M\"{o}bius function $\mu$ : $L_{\mathcal A} \longrightarrow {\bf Z}$ is
defined
by $$\begin{array}{*{3}c}
\mu(\hat{0}) & = & 1\\
\mu(t) & = & -\sum\limits_{s < t}\mu(s) \mbox{, if } \hat{0}< t
\end{array}$$
\end{defn}

\begin{defn}
The Poincar\'e polynomial $\poina$ and characteristic polynomial $\chia$
are defined by:
$$ \poina = \sum\limits_{X \in L_{\mathcal A}}\mu(X) \cdot (-t)^{rank(X)},
\mbox{    }
\chia = \sum\limits_{X \in L_{\mathcal A}}\mu(X) \cdot t^{dim(X)}$$
\end{defn}
The two polynomials are related via $\chia = t^{n+1} \cdot \pi({\cal A},
-t^{-1})$. Let $S = {\rm Sym}
(V^{*})$, $\underline{m}$ the irrelevant maximal ideal,
 $K$ the fraction field of $S$
 and suppose ${\mathcal A}$ consists of $d$
distinct
hyperplanes in $V$. For each hyperplane $H_i$ of ${\mathcal A}$, fix $l_i$
a nonzero linear form vanishing on $H_i$ and put $Q=\prod\limits_1^d l_i$.
 Denote the module of $p$ differentials over $k$ of $S$ and $K$ by
$\Omega^p_S$ and $\Omega^p_K$, respectively, and let ${\rm Der}\,_k S$
denote the
module of $k$ derivations of $S$.

\begin{defn}
$\Dp$ is the submodule of $\wps$ defined by
$$\Dp = \{ \theta \in \wps \mbox{ } | \mbox{ }\theta(Q,f_2,\ldots, f_p) \in
(Q),\, \forall f_i \in S\}.$$
$\Omega^p({\cal A})$ is the submodule of $\Omega^p_K$ defined by
$$\Omega^p({\cal A})=\{\omega\in\Omega^p_K \mbox{ } | \mbox { }
Q\,\omega\in\Omega^p_S\,\mbox{ and }\, Q\,d\omega\in\Omega^{p+1}_S\}.$$
\end{defn}

\noindent $\Do$ is usually called the module of ${\cal A}$ derivations,
while $\Omega^p({\cal A})$ is called the module of logarithmic
$p$ forms with poles along $\cal A$.
 When the arrangement is clear from the context,
we will drop it from the notation. Note that we have $D^0({\cal A})=
\Omega^0({\cal A})=S$ and we make the convention $D^p({\cal A})=
\Omega^p({\cal A})=0$, for $p<0$. 

If char $k \not\vert  d$, then $D^1({\cal A})$
 has a direct sum decomposition as $D^1_0 \oplus S(-1)$, where $D^1_0$
is the kernel of
the Jacobian matrix of $Q$ and $S(-1)$ corresponds to the Euler
derivation. Correspondingly, we have a decomposition
$\Omega^1\simeq\Omega^1_0
\oplus S(1)$.

Since all these modules are
graded, we may consider the  corresponding sheaves on ${\bf P}^n$, written
as usual as $\widetilde{\Omega^1}$ for the sheaf associated to $\Omega^1$.
An arrangement is called generic if for every $H_1,\ldots,H_m\in {\cal A}$,
with $m\leq n+1$, $\rank(H_1\cap\ldots\cap H_m)=m$. For a generic
arrangement $\cal A$, the sheaves
$\logforms$ and $\widetilde{\Omega^1}$ are related as follows:
after a change of coordinates, we may assume the first $n+1$ hyperplanes
are the coordinate hyperplanes; for $i \in \{1,\ldots,d-n-1\}$
write $l_i=\sum_{j=0}^n a_{i,j}x_j$ for the remaining hyperplanes. In \cite{z},
Ziegler gives a free resolution for $\Omega^1$ for a generic arrangement:
$$0\longrightarrow S^{d-n-1} \stackrel{\tau}{\longrightarrow} S(1)^d
\longrightarrow \Omega^1 \longrightarrow 0,$$
where $\tau$ is given by
$$\left( \begin{array}{*{4}c}
a_{1,0}x_0 &\cdots&\cdots & a_{d-n-1,0}x_0 \\
 \vdots & \ddots&\ddots & \vdots \\
a_{1,n}x_n &\cdots &\cdots& a_{d-n-1,n}x_n \\
-l_1 & 0 & \cdots & 0\\
0 & -l_2 & \cdots & 0\\
 \vdots & \ddots & \ddots & 0 \\
0 & \cdots & 0 & -l_{d-n-1}
\end{array}\right)$$

In corollary 3.4 of \cite{dk}, Dolgachev and Kapranov present $\logforms$
as the cokernel of a map $$V \otimes {\cal O}_{{\bf P}^n}(-1) \stackrel{\tau'}
{\longrightarrow} W \otimes {\cal O}_{{\bf P}^n},$$ where $V$ is the subspace
of $k^d$ consisting of relations on the linear forms defining ${\cal A}$,
$W$ is the subspace of $k^d$ orthogonal to $(1,1,\cdots, 1)$, and $\tau':
(a_1,\ldots, a_d) \longrightarrow (a_1l_1, \ldots ,a_dl_d)$. Thus, the images
of $\tau$ and $\tau'$ are isomorphic, and we have

\begin{center}
$\begin{array}{c c c c c c c c c}
0 & \longrightarrow &  S^{d-n-1} & \stackrel{\tau'}{\longrightarrow} &
S^{d-1}(1) &
\longrightarrow & \logforms(1) & \longrightarrow & 0\\

 &  & \downarrow & & \downarrow & & \downarrow & & \\
 0 & \longrightarrow &  S^{d-n-1} & \stackrel{\tau}{\longrightarrow} & S^{d}(1)
&
\longrightarrow & \Omega^1 & \longrightarrow & 0\\
 &  & \downarrow & & \downarrow & & \downarrow & & \\
 & &0 &\longrightarrow &S(1) & \longrightarrow & S(1) & \longrightarrow & 0
\end{array}$
\end{center}
By the snake lemma, we have $\logforms(1) \simeq \widetilde\Omega^1_0$.
This also follows from the local description, but the above makes explicit
the different gradings.

For every $X$ in the intersection lattice $L_{\cal A}$, the subarrangement
${\cal A}_X$ of $\cal A$ is defined by
$${\cal A}_X=\{H \in {\cal A} \mbox{ } | \mbox{ } X \subset H \}.$$ In general
this is not an essential arrangement, but we can write it as
${\cal A}_X\simeq{\cal A}'_X\times \Phi$, where ${\cal A}'_X$
is essential and $\Phi$ is an empty arrangement.

The functors on the intersection lattice $X\longrightarrow D^p({\cal A}_X)$
and $X\longrightarrow \Omega^p({\cal A}_X)$ are local (see Orlik and
Terao \cite{ot}, Chapter 4.6). What we will use is the fact that if
$\underline{q}$ is a prime ideal in $S$ and
$X=\cap_{\alpha_H\in {\underline{q}}}H$,
then we have a canonical isomorphism:
$$D^p({\cal A})_{\underline q}\simeq D^p({\cal A}_X)_{\underline q}$$
and a similar one for $\Omega^p$.

By definition, an arrangement $\cal A$ is free if $D^1({\cal A})$
is a free $S$-module. Following Yuzvinsky \cite{y2}, we will say that
an (essential, central) arrangement $\cal A$ is locally free if
for every $X\in L_{\cal A}$ with ${\rm rank}\,X<\,{\rm dim}\,V$,
the arrangement ${\cal A}_X$ is free.

 For a graded module $M$, let
$P(M,X)$ be its Hilbert series. There is a beautiful relation between the
modules
$\Dp$ and the characteristic polynomial:

\begin{thm}\label{thm:sto}(Solomon and Terao,\cite{st})
$$\chi({\cal A},t) = (-1)^{n+1} lim_{X \rightarrow 1} \sum_{p \ge 0}
P(\Dp;X)(t(X-1)-1)^p.$$
\end{thm}

There is a dual version of this theorem, which replaces $\Dp$
with $\opa$.
In certain situations, not all the modules $\Dp$ are needed to compute
$\poina$; the paradigm for this is the case of free arrangements.

\begin{thm}\label{thm:ter}(Terao,\cite{t})
\noindent If $\Do$ is free, then $\poina = \prod(1+a_it)$, where the $a_i$
are the degrees of the
generators of $\Do$.
\end{thm}

For arrangements on ${\bf P}^2$, $\widetilde\Omega^1$ is always locally free,
and suffices to determine the Poincar\'e polynomial (\cite{s}).
For every coherent sheaf $\cal F$ on ${\bf P}^n$, we denote by
$H^0_*({\cal F})$ the $S$-module $\oplus_{m\in{\bf Z}}H^0({\bf P}^n,
{\cal F}(m))$. Motivated by \thmref{thm:sto} we prove

\noindent {\bf Theorem} {\it {\mbox{ }}
For every rank $r$ bundle $\cal E$ on ${\bf P}^n$,
$$ c_t({\cal E}) = lim_{X\rightarrow 1}(-1)^r t^r
(1-X)^{n+1-r}\sum_{i=0}^rP(H^0_*(\wedge^i{\cal E}); X)\cdot({{X-1}\over
t}-1)^i.$$}
\vskip .1in
As a consequence, we get the following generalization of \thmref{thm:ter}:

\noindent {\bf Theorem} {\it {\mbox{ }}
If $\cal A$ is an arrangement such that $\widetilde{\Omega^1}$
is locally free and $\poinaa$ is the class of $\poina$ in
${\bf Z}[t]/(t^{n+1})$, then
$$\poinaa=c_t(\widetilde{\Omega^1}).$$}
\vskip .1in

We characterize those arrangements for which $\widetilde\Omega^p$ is a bundle:

\noindent {\bf Theorem}{\it {\mbox{  }}
$\widetilde\Omega^p$ is a bundle iff for every $X\in L_{\cal A}$ with ${\rm
rank}\, X<{\rm dim}\,V$,
$\Omega^p({\cal A}_X)$ is free.}
\vskip .1in
% The class of arrangements for which the conditions above hold includes
%(obviously) free arrangements and arrangements in ${\bf P}^2$, in fact, it
%is {\it much} larger,
%and  includes generic arrangements. At the end of the paper, we give an
%example of a well known
%arrangement in ${\bf P}^3$ which is neither free nor generic, but which
%is locally free.
%
We will use freely results from commutative
algebra for which our main reference is Eisenbud
\cite{e}, as well as results about Chern classes of vector bundles on
projective space, for which we refer to Fulton \cite{fu}.

\vskip .1in

\noindent{\bf Acknowledgements} We are very grateful to David Eisenbud
for useful discussions.

\section{Locally Free Arrangements}\label{sec:lfa}

We start with a general lemma about the depth of the modules
$\Omega^p$ and $D^p$.

\begin{lem}\label{lem:depth}
For every central arrangement $\cal A$ in $V$, with ${\rm dim}\,V=n+1\geq 2$,
we have ${\rm depth}\,\Omega^p\geq 2$ and ${\rm depth}\,D^p\geq 2$, for every
$p$, $1\leq p\leq n+1$.
\end{lem}

\begin{proof}
We consider the case of the modules $\Omega^p$. Recall that if $K$
is the quotient field of $S$, then
$$\Omega^p=\{\omega\in\wedge^p\Omega_S\otimes_SK\,\vert\,Q\omega\in\wedge^p
\Omega_S, Q\,d\omega\in\wedge^{p+1}\Omega_S\},$$
where $Q$ is the product of the linear forms defining the elements of
$\cal A$. In particular, $\Omega^p$ is torsion-free and therefore
${\rm depth}\,\Omega^p\geq 1$.

We have to prove that if $0\neq a\in S$ and $\omega\in\Omega^p$ are such
that $\underline{m}\,\omega\subset a\,\Omega^p$, then $\omega\in a\,\Omega^p$.
Note that since ${\rm depth}\,S\geq 2$ and the $S$-modules $\Omega^i_S$
are free, if for $\tau\in\Omega^i_K$ and $0\neq b\in S$
we have $\underline{m}\,\tau\subset b\,\Omega^i_S$, then $\tau\in\Omega^i_S$.

By definition, we have $\underline{m}\,Q\omega\subset a\,\Omega^p_S$
and the above observation gives $Q\omega/a\in\Omega^p_S$.
For every $f\in\underline{m}$ we have also $Q\,d(f\omega/a)\in\Omega^{p+1}_S$.
We use
$$Q\,d(f\omega/a)=df\wedge Q\,\omega/a+Qfd(\omega/a).$$
Since we have already seen that $Q\,\omega/a\in\Omega^p_S$ we obtain
$\underline{m}\,Qd(\omega/a)\subset\Omega^{p+1}_S$. One more application
of the above observation gives $Q\,d(\omega/a)\in\Omega^{p+1}_S$ and therefore,
$\omega/a\in\Omega^p$, which completes the proof.

The proof of the fact that ${\rm depth}\,D^p\geq 2$ is similar, using the
definition of this module.
\end{proof}

The following Proposition is the generalization of Theorem 4.75
in Orlik and Terao \cite{ot} which is the case $p=1$.
Though the general result seems known to experts, we include the proof, as
we could not find a reference in the literature.

\begin{prop}\label{prop:dual}
For every central arrangement $\cal A$, each of the 
modules $D^p$ and $\Omega^p$ is dual to the other. 
\end{prop}

\begin{proof}
A standard generalization of the argument in Orlik and Terao \cite{ot},
Proposition 4.74, gives a bilinear map of $S$-modules
$\Omega^p\times D^p\longrightarrow S$, which induces morphisms
$\alpha\,:\,\Omega^p\longrightarrow {\rm Hom}_S(D^p,S)$ and 
$\beta\,:\,D^p\longrightarrow {\rm Hom}_S(\Omega^p,S)$. The proofs of the fact
that $\alpha$ and $\beta$ are isomorphisms are similar, so that we will
give the proof only for $\alpha$. We make induction on $n\geq 1$,
the case $n=1$ being straightforward.

\lemref{lem:depth} gives ${\rm depth}\,\Omega^p\geq 2$, while it is
an easy exercise to see that since $n\geq 1$, for every graded $S$-module $M$,
${\rm depth}\,{\rm Hom}_S(M,S)\geq 2$. Since for a module $N$
of depth at least two, $N\simeq\sum_tH^0(\tilde{N}(t))$, in order
to prove that
$\alpha$ is an isomorphism, it is enough to prove that it is an
isomorphism at the sheaf level i.e. $\alpha\otimes S_{\underline q}$
is an isomorphism for every prime ideal $\underline{q}\neq\underline{m}$.
Note that if the Proposition is true for an arrangement $\cal A$
and every $p$,
then it is true also for ${\cal A}\times\Phi$, where $\Phi$
is the empty arrangement in $k$.
Indeed, if $S'={\rm Sym}(V^*\times k)$ is the polynomial ring corresponding
to ${\cal A}\times\Phi$, then by Orlik and Terao \cite{ot}, Proposition 4.84
and Solomon and Terao \cite{st}, Proposition 5.8,
 we have canonical isomorphisms 
$$\Omega^p({\cal A}\times\Phi)\simeq(\Omega^p({\cal A})\otimes_SS')\oplus
(\Omega^{p-1}({\cal A})(-1)\otimes_SS')$$
$$D^p({\cal A}\times\Phi)\simeq (D^p({\cal A})\otimes_SS')
\oplus (D^{p-1}({\cal A})(1)\otimes_SS').$$
Therefore, it follows by induction that we may assume $\cal A$
to be essential.

For a prime ideal $\underline{q}\neq\underline{m}$,
 if we take $X=\cap_{\alpha_H\in\underline{q}}H$,
then $X\in L_{\cal A}$, and ${\rm rank}\,X<{\rm dim}\,V$, since
$\cal A$ is essential. But since $\Omega^p(-)$ and $D^p(-)$
are local functors, we have canonical isomorphisms
$D^p_{\underline{q}}\simeq D^p({\cal A}_X)_{\underline{q}}$ and
$\Omega^p_{\underline{q}}\simeq \Omega^p({\cal A}_X)_{\underline{q}}$.
Since ${\cal A}_X$ is not essential, we have seen that it satisfies the
conclusion of the Proposition, and therefore we get $\alpha\otimes
S_{\underline q}$ isomorphism.
\end{proof}

\begin{thm}\label{thm:eqchar}

For an (essential, central) arrangement $\cal A$ and every positive
 integer $p$, the following are equivalent:

1. $\widetilde{D^p}$ is locally free on ${\bf P}^n$.

1'. $\widetilde{\Omega ^p}$ is locally free on ${\bf P}^n$.

2. For every $X\in L_{\cal A}$ with ${\rm rank}\,X<{\rm dim}\,V$,
$D^p({\cal A}_X)$ is free.

2'. For every $X\in L_{\cal A}$ with ${\rm rank}\, X<{\rm dim}\,V$,
$\Omega^p({\cal A}_X)$ is free.
\end{thm}

\begin{proof} The equivalences
$1\Leftrightarrow1'$
 and $2\Leftrightarrow2'$ follow from \propref{prop:dual}.

For the proof of $1\Rightarrow2$,
 let $X\in L_{\cal A}$ be a nonzero linear subspace and
$I_X\subset S$ its ideal. By making a linear change of variables we can
assume that $I_X=(X_0,\ldots,X_{r-1})$, where $r={\rm rank}\,X$.
Since ${\cal A}_X\simeq {\cal A'}_X\times\Phi^{n+1-r}$, if
$S_1=k[X_0,\ldots,X_{r-1}]$, then
$$D({\cal A}_X)\simeq(D({\cal A'}_X)\otimes_{S_1}S)\oplus S^{n+1-r}$$
and more generally
$$D^p({\cal A}_X)\simeq\bigoplus_{0\leq i\leq p}(D^i({\cal A'}_X)\otimes
S^{{n+1-r}\choose {p-i}}).$$

Since $D^i({\cal A'}_X)$ is a free $S_1$ module if and only if
$D^i({\cal A'}_X)_{I_X}$ is a free $S_1$-module, it follows that
$D^p({\cal A}_X)$ is free if and only if $D^p({\cal A}_X)_{I_X}$
is free. But because $\widetilde{D^p}$ is locally free and
$I_X\neq \underline{m}$, $D^p_{I_X}$ is free over $S_{I_X}$.
On the other hand, since $D^p(-)$ is a local functor we have
$$D^p_{I_X}\simeq D^p({\cal A}_X)_{I_X},$$
and therefore $D^p({\cal A}_X)$ is free.

In order to prove $2\Rightarrow1$, let us consider a prime ideal
$\underline{q}$ different from $\underline{m}$. If we take
$X=\cap_{\alpha_H\in\underline{q}}H$, then
$X\in L_{\cal A}$, ${\rm rank}\, X<{\rm dim}\, V$ and because
$D^p(-)$ is a local functor we have $D^p_{\underline{q}}
\simeq D^p({\cal A}_X)_{\underline{q}}$, which is free over $S_{\underline{q}}$
by hypothesis. This concludes the proof of the theorem.
\end{proof}

By taking $p=1$ in the above theorem we obtain the following:

\begin{cor}\label{charlf}
An arrangement $\cal A$ is locally free if and only if $\widetilde{D^1}$
is locally free on ${\bf P}^n$.
\end{cor}

\begin{rem}\label{terao}
A famous conjecture due to Terao asserts that the freeness of an arrangement
depends only on the intersection lattice. This is equivalent to the fact
that the local freeness of an arrangement depends only on
the intersection lattice. Indeed, the fact that the second statement is a
consequence of
Terao's conjecture follows immediately by induction on
rank. Conversely, if two arrangements ${\cal A}_1$ and ${\cal A}_2$
have isomorphic lattices and ${\cal A}_1$ is free, consider the
product arrangements ${\cal A}'_1={\cal A}_1\times B_1$
 and ${\cal A}'_2={\cal A}_2\times B_1$,
where $B_1$ is the Boolean arrangement in $W=k$. Then we have
${\cal A}'_1$ free and in particular locally free,
 while ${\cal A}'_2$
is free if and only if it is locally free if and only if
${\cal A}_2$ is free.

Notice also that since free arrangements are always locally free,
Terao's conjecture becomes a question on the splitting of vector bundles
on ${\bf P}^n$.
\end{rem}

\begin{lem}\label{aa'}
For every arrangement $\cal A$, every $p\geq 1$ and every $X\in L_{\cal A}$,
$D^p({\cal A}_X)$ is free if and only if $D^i({\cal A}'_X)$ is free,
for all $i$ with $p-{\rm dim}\, X\leq i\leq p$.
\end{lem}

\begin{proof}
We have seen in the proof of the above \thmref{thm:eqchar} that we can write
$$
D^p({\cal A}_X)=\bigoplus_{0\leq i\leq p}
 D^i({\cal A}'_X)\otimes_{S_1}S^{{n+1-r}
\choose {p-i}},
$$
where $S_1=k[X_0,\ldots,X_{r-1}]$ and $r=n+1-{\rm dim}\, X$.
To conclude it is enough to notice that
$S^{{n+1-r}\choose {p-i}}\neq 0$ if and only if
$p-{\rm dim}\, X\leq i\leq p$.
\end{proof}

\begin{cor}
For every arrangement $\cal A$ and every $X\in L_{\cal A}$ such that
${\rm dim}\, X\geq p-1$, $D^p({\cal A}_X)$ is free if and only if
${\cal A}_X$ is free. In particular, if $\widetilde{D^p}$ is locally free,
then for every $X\in L_{\cal A}$ with ${\rm dim}\, X\geq \max\{1,p-1\}$,
${\cal A}_X$ is free.
\end{cor}

\begin{proof}
The first assertion follows from the above lemma, since we have
$p-{\rm dim}\,X\leq 1$. The second assertion follows from the first one and
\thmref{thm:eqchar}.
\end{proof}

\begin{cor}
For every arrangement $\cal A$, $\widetilde{D^1}$ is locally free
if and only if $\widetilde{D^2}$ is locally free.
\end{cor}

\begin{proof}
The ``if'' part follows from the previous corollary in the case $p=2$,
while the ``only if'' part is a consequence of the more general
proposition below.
\end{proof}

\begin{prop}\label{prop:islocfree}
If $\widetilde{D^1}$ is locally free, then the natural map
$\wedge^p(D^1)\longrightarrow D^p$ induces an isomorphism
$\wedge^p(\widetilde{D^1})\simeq\widetilde{D^p}$ for every $p\geq 1$.
In particular, $\widetilde{D^p}$ is locally free. The similar
assertion about the natural morphism $\wedge^p(\Omega^1)\longrightarrow
\Omega^p$ is also true.
\end{prop}

\begin{proof}
We have to prove that for every $\underline{q}\in{\rm Spec}\,(S)$,
$\underline{q}\neq\underline{m}$, the localized morphism:

$$
\wedge^p(D^1_{\underline{q}})\longrightarrow D^p_{\underline q}
$$
is an isomorphism. Since $D^1(-)$ and $D^p(-)$ are local functors,
if $X=\cap_{\alpha_H\in\underline{q}}H$, we have
$D^1_{\underline{q}}=D^1({\cal A}_X)_{\underline{q}}$ and
$D^p_{\underline{q}}=D^p({\cal A}_X)_{\underline{q}}$.

Because $\underline{q}\neq\underline{m}$, $X\neq (0)$ and by
hypothesis ${\cal A}_X$ is free. But for free arrangements
$D^p\simeq\wedge^p(D^1)$ (see Solomon-Terao \cite{st}, Prop.3.4),
which concludes the proof of the first assertion.

The proof of the last statement is similar.
\end{proof}

When ${\rm char}\,k$ does not divide $|{\cal A}|$, it is possible to
characterize the local freeness of $\cal A$ using the Ext
modules of the Jacobian ideal.

\begin{prop}\label{prop:Ext}
For an arrangement $\cal A$ such that ${\rm char}\,k\not\vert |{\cal A}|$,
$\cal A$ is locally free if and only if the modules
${\rm Ext}^i_S(S/J,S)$ are supported only at the maximal ideal
$\underline m$, for all $i\geq 3$.
\end{prop}

\begin{proof}
Recall that we have the module $D_0$ defined by the exact sequence:
$$
0\longrightarrow D_0\longrightarrow S^{n+1}\longrightarrow
 J(d-1)\longrightarrow 0,
$$
where $d=|{\cal A}|$ and $J$ is the Jacobian ideal.
Since ${\rm char}\,k\not\vert d$, we have
$D\simeq D_0\oplus S(-1)$.
Using \thmref{thm:eqchar} it follows that $\cal A$ is locally free
if and only if $\widetilde{D_0}$ is locally free.

For a finitely generated $S$-module $M$ it is known that the set
$$
S(M)=\{\underline{q}\in{\rm Spec}(S)\,\vert\,M_{\underline{q}}\,
{\rm is}\, {\rm not}\, {\rm a}\,{\rm free}\,S_{\underline{q}}-{\rm module}\}
$$
can be written as
$$
S(M)=\bigcup_{i\geq 1}{\rm Supp}\,{\rm Ext}^i_S(M,S)
$$
(see,  Hartshorne \cite{h}, p. 238, exercise 6.6).

Therefore $\widetilde{D_0}$ is locally free if and only if
${\rm Supp}\,{\rm Ext}^i_S(D_0,S)\subset\{\underline{m}\}$
for every $i\geq 1$. Since ${\rm Ext}^i_S(S/J, S)\simeq {\rm Ext}^{i-2}
_S(D_0, S)$ for every $i\geq 3$, the proof of the proposition is complete.
\end{proof}

\section{Chern classes of vector bundles on ${\bf P}^n$}
\label{sec:vecbundles}

We consider a vector bundle $\cal E$ of rank $r$ on ${\bf P}^n$.
Motivated by the application in the context of arrangements which
will be given in the next section, we introduce
$R({\cal E};t,X)\in{\bf Z}[t]((X))/(t^{n+1})$, defined by
$$
R({\cal E};t,X)=(-1)^r t^r (1-X)^{n+1-r}\sum_{i=0}^rP(H^0_*(\wedge^i{\cal E});
X)\cdot({{X-1}\over t}-1)^i.
$$
Here $H^0_*(\wedge^i{\cal E})$ is the finitely generated graded
module $\oplus_{m\in{\bf Z}}H^0({\bf P}^n,\wedge^i{\cal E}(m))$.
Recall that for a finitely generated graded $S$-module $M$,
$P(M;X)$ denotes the Hilbert series of $M$ which is a Laurent series,
but also a rational function in $X$. The main result of this
section is that $R({\cal E};t,X)$
can be used to compute the Chern polynomial of $\cal E$. More
precisely, we have the following

\begin{thm}\label{thm:limR}
If $\cal E$ is a vector bundle on ${\bf P}^n$, then
$$
\lim_{X\to 1} R({\cal E};t, X)=c_t({\cal E}).
$$
\end{thm}

\begin{rem}\label{rem:indep}
We will see in the proof that
 in order to compute
the above limit, in the definition of $R({\cal E};t,x)$ we may replace
each $H^0_*(\wedge^i{\cal E})$ with a different finitely generated
module $M_i$ such that $\wedge^i{\cal E}\simeq\widetilde{M_i}$.
\end{rem}

Before proving the theorem we give two lemmas.

\begin{lem}\label{lem:rk1}
The assertion of \thmref{thm:limR} is true in the case of a split
vector bundle $\cal E$.
\end{lem}

\begin{proof}
Suppose that ${\cal E}\simeq{\cal O}(a_1)\oplus\ldots\oplus{\cal O}(a_r)$.
In this case we have
$$
R({\cal E};t,X)=(-1)^rt^r(1-X)^{n+1-r}\sum_{i=0}^r
P(H^0_*(\oplus_{1\leq k_1<\ldots<k_i\leq r}{\cal O}(a_{k_1}+\ldots+a_{k_r});
X)\cdot({{X-1}\over t}-1)^i.
$$
Since for every $a\in{\bf Z}$, $H^0_*({\cal O}(a))=S(a)$ and
$P(S(a); X)=X^{-a}(1-X)^{-n-1}$, we get
$$
R({\cal E};t,X)=(-1)^rt^r(1-X)^{n+1-r}\sum_{i=0}^r\sum
_{1\leq k_1<\ldots<k_i\leq r}X^{-a_{k_1}-\ldots-a_{k_i}}
\cdot(1-X)^{-n-1}\cdot({{X-1}\over t}-1)^i
$$
$$=(-1)^rt^r(1-X)^{-r}\prod_{i=1}^r(1+X^{-a_i}({{X-1}\over t}-1))
=(-1)^r\prod_{i=1}^r(t\cdot{{1-X^{-a_i}}\over {1-X}}-X^{-a_i}).
$$
It follows from this that the limit exists and
$$
\lim_{X\to 1}R({\cal E};t,X)=(-1)^r\prod_{i=1}^r(-a_it-1)
=\prod_{i=1}^r(1+a_it)=c_t({\cal E}).
$$
\end{proof}

\begin{lem}\label{lem:polcc}
If $r\geq n$ is fixed and $P\in{\bf Q}[X_1,\ldots,X_n]$ is a polynomial
such that
$$P(c_1({\cal E}),\ldots,c_n({\cal E}))=0,$$
for every split vector bundle of rank $r$, ${\cal E}={\cal O}(a_1)
\oplus\ldots\oplus{\cal O}(a_r)$, then $P=0$.
\end{lem}

\begin{proof}
If $s_i$, $1\leq i\leq n$ is the $i^{\rm th}$ symmetric polynomial
in $a_1,\ldots,a_r$, then $c_i({\cal E})=s_i(a_1,\ldots,a_r)$, for
every split vector bundle $\cal E$ as in the hypothesis. Therefore
we get $P(s_1,\ldots,s_n)=0$.

But the morphism
$$\phi\,:\,{\bf Q}[X_1,\ldots,X_n]\longrightarrow {\bf Q}[Y_1,\ldots,Y_r]
$$
given by $\phi (X_i)=s_i(Y_1,\ldots,Y_r)$ is the restriction to
${\bf Q}[X_1,\ldots,X_n]$ of the monomorphism
$$\psi\,:\,{\bf Q}[X_1,\ldots,X_r]\longrightarrow {\bf Q}[Y_1,\ldots,Y_r],$$
where $\psi(X_i)=s_i(Y)$, for $1\leq i\leq r$. Therefore we have $P=0$.
\end{proof}

We are now ready to give the proof of the theorem.

\begin{proof}
It is easy to check that if ${\cal E}'={\cal E}\oplus{\cal O}$,
then $R({\cal E}';t,X)=R({\cal E};t,X)$. As we have also
$c_t({\cal E}')=c_t({\cal E})$, it follows that by taking the direct sum with
a large enough number of trivial bundles, we may suppose that ${\rm rank}\,
{\cal E}=r\geq n$.

In this case, using the above two lemmas, we see that in order to prove
the theorem it is enough to show that the existence of the limit
and its value
 can be expressed
in terms of some
polynomial identities with rational coefficients
in the Chern classes of $\cal E$. We have
$$
R({\cal E};t,X)=(-1)^r
t^r (1-X)^{n+1-r}\sum_{i=0}^r P(H^0_*(\wedge^i{\cal E});
X)\cdot\sum_{j=0}^i(-1)^{i-j}{i\choose j}({{X-1}\over t})^j
$$
$$
=(-1)^r\sum_{j=0}^rt^{r-j}(1-X)^{n+1-r+j}\cdot\sum_{i=j}^r(-1)^i{i\choose
j}P(H^0_*(\wedge^i{\cal E}); X).
$$

Since for $0\leq i\leq r$, $H^0_*(\wedge^i{\cal E})$ is a $S$-module
of dimension $n+1$, we can write
$$
P(H^0_*(\wedge^i({\cal E}));X)={{Q_i(X)}\over{(1-X)^{n+1}}},
$$
where $Q\in{\bf Z}[X,X^{-1}]$.
Moreover, if we consider the Taylor expansion of $Q_i(X)$ around $X=1$:
$$
Q_i(X)=\sum_{l\geq 0}e^{(i)}_l(X-1)^l,
$$
then the first $n+1$ coefficients of this expansion can be recovered
from the Hilbert polynomial of $H^0_*(\wedge^i{\cal E})$, which can
be written as:
$$
T_i(X)=\sum_{l=0}^n(-1)^{n-l}e^{(i)}_{n-l}{{X+l}\choose l}.
$$
For these results, see for example Bruns and Herzog \cite{bh}, Chapter 4.1.

Therefore  we have
$e^{(i)}_k=(-1)^k\Delta^{n-k}\,T_i(0)$. For every graded $S$-module
$M$, its  Hilbert polynomial is given by the formula
$T(m)=\chi(\widetilde{M}(m))$, for every $m\in{\bf Z}$. Using the
short exact sequences corresponding to successive hyperplane sections
we get
$$e^{(i)}_k=(-1)^k\chi(\wedge^i{\cal E}\vert_{H_k}),$$
where $H_k\subset {\bf P}^n$ is a linear subspace of dimension $k$,
for $0\leq k\leq n$. We deduce

$$
R({\cal E};t,X)=(-1)^r\sum_{j=r-n}^rt^{r-j}(1-X)^{n+1-r+j}\cdot
\sum_{i=j}^r(-1)^i{i\choose j}\cdot{{\sum_{k\geq 0}(-1)^ke^{(i)}_k(1-X)^k}
\over {(1-X)^{n+1}}}.
$$
$$
R({\cal E};t,X)=(-1)^r\sum_{j=r-n}^rt^{r-j}\cdot\sum_{i=j}^r
(-1)^i{i\choose j}\sum_{k\geq 0}(-1)^k{{e^{(i)}_k}\over {(1-X)^{r-j-k}}}.
$$

By considering the coefficient of $t^{r-j}$, for $r-n\leq j\leq r$,
the fact that $\lim_{X\to 1}R({\cal E};t,X)$ exists and is
equal to $c_t({\cal E})$ is equivalent to:

$$(1)\,\,\sum_{i=j}^r(-1)^i{i\choose j}e^{(i)}_k=0,\,{\rm for}\,0\leq k\leq
r-j-1,
$$
$$
(2)\,\,\sum_{i=j}^r(-1)^i{i\choose j}e^{(i)}_{r-j}=(-1)^rc_{r-j}({\cal E}),
$$
for every $j$ with $r-n\leq j\leq r$.

Using the formulas we have for $e^{(i)}_k$, these relations become:
$$
(1')\,\,\sum_{i=j}^r(-1)^i{i\choose j}\chi(\wedge^i{\cal
E}\vert_{H_k})=0,\, {\rm for}\,0\leq k\leq r-j-1
$$
$$
(2')\sum_{i=j}^r(-1)^i{i\choose j}\chi(\wedge^i{\cal E}\vert_{H_{r-j}})=
(-1)^rc_{r-j}({\cal E}),
$$
for every $j$, with $r-n\leq j\leq r$.

For future reference, notice that for $j=r-n$, $(2')$ becomes
$$
(2'')\,\,\sum_{i=r-n}^r(-1)^i{i\choose {r-n}}\chi(\wedge^i{\cal E})
=(-1)^rc_n({\cal E}).
$$

In order to finish the proof of the theorem, it is enough to notice
that using Hirzebruch-Riemann-Roch theorem (see Fulton \cite{fu}
Corollary 15.2.1) all the Euler-Poincar\'e characteristics can
be expressed as polynomials with rational coefficients in the
Chern classes of the exterior powers $\wedge^i{\cal E}$. Indeed,
since the Chern classes of these exterior powers can be computed
as polynomials in the Chern classes of $\cal E$, we can apply
\lemref{lem:polcc} and \lemref{lem:rk1} to conclude the proof
of the theorem.
\end{proof}

\begin{cor}\label{cor:ccid}
If $\cal E$ is a vector bundle on ${\bf P}^n$ with ${\rank}\,{\cal E}=r
\geq n$, then we have the following formula for the top Chern class:
$$\sum_{i=r-n}^r(-1)^i{i\choose {r-n}}\chi(\wedge^i{\cal E})=(-1)^rc_n
({\cal E}).$$
In particular, if $r=n$ we have
$$\sum_{i=0}^n(-1)^i\chi(\wedge^i{\cal E})=(-1)^nc_n({\cal E}).$$
\end{cor}

\begin{proof}
This is just the identity $(2'')$ in the proof of \thmref{thm:limR}.
\end{proof}

\begin{rem}\label{lem:degenlocus}
If $\cal E$ is a vector bundle on ${\bf P}^n$ such that ${\cal E}^*$
has $r-n+1$ sections $\sigma_1,\ldots,\sigma_{r-n+1}$ such that the
degeneration locus $\Gamma$ is zero dimensional and the degeneration
locus $\Gamma'$ of $\sigma_1,\ldots,\sigma_{r-n}$
 is empty, then the formula in \corref{cor:ccid}
is equivalent to the formula giving the degree of $\Gamma$.

Indeed, $\sigma_1,\ldots, \sigma_{r-n+1}$ define a morphism
${\cal E}\longrightarrow{\cal F}={\cal O}_{{\bf P}^n}^{r-n+1}$
which gives an Eagon-Northcott type complex
$$0\longrightarrow\wedge^r{\cal E}\otimes(S_n{\cal F})^*
\longrightarrow\ldots\longrightarrow\wedge^{r-n+1}{\cal E}\otimes{\cal F}^*
\longrightarrow\wedge^{r-n}{\cal E}.$$

It follows from Eisenbud \cite{e}, Theorems A.2.10 and A.2.14 that
since ${\dim}\,\Gamma=0$ and $\Gamma'=\emptyset$, this complex is exact
and moreover, it gives a resolution of ${\cal O}_{\Gamma}$.
Therefore we obtain $$\sum_{i=r-n}^r(-1)^i\chi(\wedge^i{\cal E})=
(-1)^{r-n}\chi({\cal O}_{\Gamma})=(-1)^{r-n}{\rm deg}\,\Gamma.$$

On the other hand, the Thom-Porteous formula (see Fulton \cite{fu},
Theorem 14.4) says that under our hypothesis
${\rm deg}\,\Gamma=(-1)^nc_n({\cal E})$ and we get the formula in
\corref{cor:ccid}.
\end{rem}

Let ${\cal E}$ be a vector bundle on ${\bf P}^n$, with ${\rm rank}\,{\cal E}
=n$. For every $i$, we denote by $Q(\wedge^i{\cal E}; X)$ the Hilbert
polynomial of $\wedge^i{\cal E}$.

\begin{cor}
With the above notation, we have
$$\sum_{i=0}^n(-1)^iQ(\wedge^i({\cal E}); iX)=(-1)^n\sum_{i=0}^nc_i({\cal E})
X^{n-i}.$$
\end{cor}

\begin{proof}
We prove that the two polynomials take the same value for every $a\in\Z$.
Indeed, if in \corref{cor:ccid} we replace ${\cal E}$
by ${\cal E}(a)$, then we have
$$\chi(\wedge^i({\cal E}(a)))=\chi(\wedge^i{\cal E})(ai)=
Q(\wedge^i{\cal E}; ai),$$
while $c_n({\cal E}(a))=\sum_{i=0}^nc_i({\cal E})a^{n-i}$ (see Fulton
\cite{fu}, Remark 3.2.3).
\end{proof}

\section{The Characteristic Polynomial}\label{sec:charpoly}

In this section we will apply the general results we have obtained so far
to the case of the vector bundle associated to the module of
$\cal A$ derivations of a locally free arrangement ${\cal A}$.
Recall that $\poina$ is a polynomial of degree $n+1$. We will
denote by $\poinaa$ its class in $\Z[t]/(t^{n+1})$.
The main result is the following.

\begin{thm}\label{thm:char}
If $\cal A$ is a locally free arrangement, then
$$\poinaa=c_t(\widetilde{\Omega^1}).$$
\end{thm}

\begin{proof}
From the basic relation between the Poincar\'e and the characteristic
polynomial we get
$$\poina=(-t)^{n+1}\chi({\cal A},-t^{-1}).$$
Combining this with \thmref{thm:sto} we obtain
$$\poina=\lim_{X\to 1}t^{n+1}\sum_{i\geq 0}
P(D^i;X)(-{1\over t}(X-1)-1)^i.$$
By \propref{prop:islocfree} we have
$\widetilde{D^i}=\wedge^i\widetilde{D^1}$;
so from \remref{rem:indep} it follows that
$$\lim_{X\to 1}R(\widetilde{D^1};-t, X)=\lim_{X\to 1}
\sum_{i=0}^nP(D^i; X)\cdot(-{1\over t}(X-1)-1).$$
By \thmref{thm:limR} this limit is equal to
 $$c_{-t}(\widetilde{D^1})=c_t(\widetilde{\Omega^1}).$$
We therefore obtain $$\poinaa=c_t(\widetilde{\Omega^1}).$$
\end{proof}

\begin{rem}
Since it is known that $\pi({\cal A};-1)=0$ (see Orlik and Terao \cite{ot},
Proposition 2.5.1), in order to know $\poina$ it is enough to know $\poinaa$.
\end{rem}

In fact, when ${\rm char}\,k$ does not divide $|{\cal A}|$, then
from \thmref{thm:char}
we can deduce a formula for $\poina$ involving the vector
bundle $\widetilde{\Omega^1_0}$.

\begin{cor}\label{cor:poinchern}
If ${\cal A}$ is a locally free arrangement such that ${\rm char}\,k$
does not divide $|{\cal A}|$, then we have
$$\poina=(1+t)c_t(\widetilde{\Omega^1_0}).$$
\end{cor}

\begin{proof}
Since $\poina /(1+t)$ is a polynomial of degree $n$, it follows that
it is enough to prove that its class
in $\Z[t]/(t^{n+1})$ is equal to $c_t(\widetilde{\Omega^1_0})$.
But by \thmref{thm:char} it follows that
this class is equal to
$$\poinaa/(1+t)=c_t(\widetilde{\Omega^1})/(1+t)=c_t(\widetilde{\Omega^1_0}),$$
since $\widetilde{\Omega^1}\simeq\widetilde{\Omega^1_0}\oplus
{\cal O}(1)$.
\end{proof}
In \cite{y2}, Yuzvinsky proves that for a locally free arrangement,
the Hilbert polynomial of the module $D^1({\cal A})$ depends
only on the lattice. The following corollary makes this more precise.

\begin{cor}\label{cor:latdep}
If ${\cal A}$ is a locally free arrangement, then giving the Hilbert
polynomial of $D^1({\cal A})$ is equivalent to giving the Poincar\'e
polynomial $\poina$ of the arrangement.
\end{cor}

\begin{proof}
The statement follows from \thmref{thm:char} and the
Hirzebruch-Riemann-Roch theorem
(see Fulton \cite{fu}, Corollary 15.2.1).
\end{proof}

\begin{exm} Let ${\cal A}$ be the arrangement in ${\bf P}^3$ defined by the
vanishing of the fifteen linear forms
$a_0x_0+a_1x_1+a_2x_2+a_3x_3$, where $a_i$ is either zero or one.
This example was constructed by Edelman-Reiner (\cite{er}) as a
counterexample to a conjecture of Orlik; $\poina = 1+15t+80t^2+170t^3+104t^4$.
There are 45 rank three elements of $L_{\cal A}$; we consider the
corresponding subarrangements as essential arrangements in ${\bf P}^2$.
The subarrangements are of three distinct types, described below
($\mu(L_2)$ is the M\"obius function of the rank two elements of $L_{{\cal
A}_X}$)
:

$\begin{array}{*{3}c}
& & \\
20 & \mbox{subarrangements on 3 hyperplanes, } & \mu(L_2) = (1,1,1)\\
& & \\
15 & \mbox{subarrangements on 5 hyperplanes, }  & \mu(L_2) = (2,2,1,1,1,1)\\
& & \\
10 & \mbox{subarrangements on 7 hyperplanes, } & \mu(L_2) =
(2,2,2,2,2,2,1,1,1)\\
& &
\end{array}$

It is easy to check that all of these subarrangements are free,
so ${\cal A}$ is locally free. We illustrate \corref{cor:latdep} for this
example. For a rank three bundle ${\cal E}$ on ${\bf P}^3$ we have:
$$\int \mbox{ch}({\cal E}(m)) \cdot \mbox{td}({\cal T}_{{\bf P}^3}) =
\frac{1}{2}m^3+(3+\frac{c_1}{2})m^2+(\frac{11}{2}+2c_1+\frac{c_1^2}{2}-c_2)m
+3+\frac{11c_1}{6}+c_1^2-2c_2+\frac{c_1^3}{3}-\frac{c_1c_2}{2}+\frac{c_3}{2}.$$
Since we know $\poina$, we may apply \corref{cor:poinchern} to obtain
the Chern classes $c_i$, and from Hirzebruch-Riemann-Roch we obtain the
Hilbert polynomial:
$$\chi(D^1_0(m)) = \frac{1}{2}m^3-4m^2+\frac{57}{6}m-6.$$ The point is that
for a locally free arrangement, knowing the combinatorial data (i.e. Poincar\'e
polynomial) means knowing the Hilbert polynomial, which can make the
computation
of the free resolution much faster; for this example the free resolution is:

$$0 \longrightarrow S(-6) \longrightarrow  S^4(-5) \longrightarrow  D^1_0
  \longrightarrow  0.$$
Thus, we see that $\widetilde{\oo} \simeq \Omega_{\bf
P^3}(6)$,  and $c_t(\widetilde{\oo}) = \frac{(1+5t)^4}{(1+6t)}
\mbox{ mod }t^4 = 1+14t+66t^2+104t^3$, as expected.
\end{exm}
{\bf Remark} In \cite{dk}, Dolgachev and Kapranov point out that for a generic
arrangement,
$\Omega^1(\mbox{log } D)$ is a  Steiner bundle (hence stable); in the
previous example
$\widetilde{D^1_0}$ is Steiner. In ${{\bf P}^2}$ there are many examples of
non-generic arrangements for which $\widetilde{\Omega^1_0}$ is
indecomposable but
not semistable.

\section{Minimal free resolutions for the modules of logarithmic forms}

In this section we give a minimal free
resolution for the modules $\Omega^p({\cal A})$
of logarithmic forms in the case of a locally free arrangement ${\cal A}$
with ${\rm pdim}\,\Omega^1({\cal A})=1$. This generalizes the results
of Rose and Terao \cite{rt} in the case of generic arrangements.
The same idea can be used to give a minimal free resolution for the modules
$D^p({\cal A})$ when ${\rm pdim}\,D^1({\cal A})=1$.

 We first recall the definition of syzygy modules:

\begin{defn}
A module $M$ is a ${\it k^{th} syzygy}$ if there exists an exact sequence
$$0\longrightarrow M \longrightarrow F_1 \longrightarrow F_2 \longrightarrow
\ldots \longrightarrow F_k,$$
with $F_i$ free.
\end{defn}

\begin{lem}\label{lem:pdimsyz}
If ${\rm pdim}\,M=1$ and $Ext^1(M,S)_{\underline q} = 0$
for every prime ideal ${\underline q} \neq {\underline m}$,
then $M$ is an $(n-1)^{st}$ syzygy.
\end{lem}
\begin{proof}

From the short exact sequence:
$$ 0 \longrightarrow F_1 \longrightarrow F_0 \longrightarrow M
\longrightarrow 0,$$
we obtain an exact sequence:
$$ 0 \longrightarrow M^* \longrightarrow F_0^* \longrightarrow F_1^*
\longrightarrow Ext^1(M,S) \longrightarrow 0.$$
Since $Ext^1(M,S)$ is supported only at $m$, $Ext^1(M,S)$ is a module of finite
length, so has a free resolution of length $n+1$; hence $M^*$ has projective
dimension $n-1$. Dualizing once more and using the fact that
$Ext^i(Ext^1(M,S),S)$
is zero for $i \ne n+1$, we find that $M$ is an $(n-1)^{st}$ syzygy.
\end{proof}

Suppose now that $M$ is a finitely generated $S$-module. Note that
if $\widetilde{M}$ is locally free, then the second condition in
\lemref{lem:pdimsyz} is automatically satisfied.
In \cite{l1} Lebelt shows that if
${\rm pdim}\, M=i$ and if $M$ is an $i(p-1)^{st}$ syzygy, then one
can obtain a free resolution for $\Lambda^pM$ from a free resolution
of $M$.
In the situation considered above, i.e. $i=1$
 and $M$ is an $(n-1)^{st}$
syzygy, we obtain a minimal free resolution of $\Lambda^p(M)$, for $p\leq n-1$,
which is an Eagon-Northcott type complex.
Namely, if
$$0\longrightarrow F_1\longrightarrow F_0\longrightarrow M\longrightarrow 0$$
is a minimal free resolution of $M$, then a minimal free resolution
of $\wedge^pM$ is given by:
$$(F_{\bullet}^{(p)})(M):\quad
 0 \longrightarrow D_pF_1 \longrightarrow D_{p-1}F_1 \otimes F_0
\longrightarrow
D_{p-2}F_1 \otimes \Lambda^2 F_0 \longrightarrow \ldots \longrightarrow
\Lambda^p F_0
\longrightarrow \Lambda^p M \longrightarrow 0,$$
where $D_iF_1=(S_i(F_1^*))^*$ denotes
 the $i^{th}$ divided power of $F_1$. If ${\rm char}\,k=0$, then
$D_iF_1\simeq S_iF_1$ is the usual symmetric power of $F_1$.

We can now give the main result of this section:

\begin{thm}\label{thm:minres}
If ${\cal A}$ is a locally free arrangement and ${\rm pdim}\,\Omega^1=1$,
then the natural morphism
$$\wedge^p\Omega^1\longrightarrow\Omega^p$$
is an isomorphism, and $(F_{\bullet}^{(p)}(\Omega^1))$
gives a minimal free resolution of $\Omega^p$, for every $p$,
$p\leq n-1$.
\end{thm}

\begin{proof}
From \lemref{lem:pdimsyz} and Lebelt's result cited above, it follows
that for every $p$, $p\leq n-1$, $F_{\bullet}^{(p)}(\Omega^1)$
is a (minimal) free resolution of $\wedge^p\Omega^1$.
In particular, we have ${\rm pdim}\,\wedge^p\Omega^1=p$.
By the Auslander-Buchsbaum formula we obtain ${\rm depth}\,\wedge^p\Omega^1=
n+1-p\geq 2$.

We consider the commutative diagram:

\begin{center}
$\begin{array}{c c c}
\wedge^p\Omega^1&\stackrel{\alpha}{\longrightarrow}&\Omega^p\\

\downarrow& &\downarrow\\

H^0_*(\wedge^p\widetilde{\Omega^1})&\stackrel{\beta}{\longrightarrow}&
H^0_*(\widetilde{\Omega^p})
\end{array}$
\end{center}

\propref{prop:islocfree} implies that $\beta$ is an isomorphism.
 By \lemref{lem:depth},
 ${\rm depth}\,\Omega^p\geq 2$ and we also have ${\rm
depth}\,\wedge^p\Omega^1
\geq 2$ and therefore both the vertical maps in the diagram are isomorphisms.
We conclude that $\alpha$ is an isomorphism.
\end{proof}

\begin{cor}\label{cor:pdimlogp}
If $\cal A$ is a locally free arrangement with ${\rm pdim}\,\Omega^1=1$,
then we have
$${\rm pdim}\,\Omega^p=p,\,{\rm for}\,p\leq n-1,$$
$${\rm pdim}\,\Omega^{n+1}=0.$$
Moreover, if ${\rm char}\,k$ does not divide $|{\cal A}|$,
then we
have also
$${\rm pdim}\,\Omega^n=n-1,$$
and therefore ${\rm pdim}\,D^1=n-1$.
\end{cor}

\begin{proof}
The first assertion follows from \thmref{thm:minres}, while the second one
is true for an arbitrary arrangement (see Orlik and Terao \cite{ot},
Proposition 4.68).

If ${\rm char}\,k\not ||{\cal A}|$, then we have
$\Omega^1=\Omega^1_0\oplus S(1)$. Therefore,
$$\wedge^p\Omega^1=\wedge^p\Omega^1_0\oplus\wedge^{p-1}\Omega^1_0(1),$$
for every $p$. It follows immediately that ${\rm pdim}\wedge^p\Omega_0=p$,
for $p\leq n-1$. In particular, from the Auslander-Buchsbaum formula
we get ${\rm depth}\,\wedge^{n-1}\Omega^1_0(1)=2$.

From this and the decomposition
 $\wedge^n\Omega^1=\wedge^n\Omega^1_0\oplus\wedge^{n-1}
\Omega^1_0(1)$ it follows that $\wedge^{n-1}\Omega^1_0(1)$ is a direct summand
of $H^0_*(\widetilde{\wedge^n\Omega})$.

By \lemref{lem:depth} we have ${\rm depth}\,\Omega^n\geq 2$.
 This implies that $\Omega^n\simeq
H^0_*(\widetilde{\Omega^n})\simeq H^0_*(\wedge^n\widetilde\Omega^1)$.
From the fact that this module has a summand of depth two and has
depth at least two, we conclude that the depth is exactly two, and the
result follows from one more application of the Auslander-Buchsbaum formula.

The fact that ${\rm pdim}\,D^1=n-1$ is a consequence of the general
fact that $\Omega^n\simeq D^1(d)$, where $d=|{\cal A}|$ (see
Rose and Terao \cite{rt}, Lemma 4.4.1).
\end{proof}

A generic arrangement $\cal A$ is locally free since for every $X\in\
L_{\cal A}$, with ${\rm rank}\,X<{\rm dim}\,V$, ${\cal A}_X$ is
isomorphic to the product between a Boolean arrangement and
an empty arrangement.
On the other hand, Ziegler \cite{z}, Corollary 7.7 shows that in this
case ${\rm pdim}\,\Omega^1=1$ and therefore \thmref{thm:minres}
and \corref{cor:pdimlogp} apply in this case and give the results
in Rose and Terao \cite{rt}.

Analogous results with similar proofs
hold if we replace the modules $\Omega^p$ with
the modules $D^p$. Namely, we have

\begin{thm}\label{thm:minres'}
If ${\cal A}$ is a locally free arrangement and ${\rm pdim}\,D^1=1$, then the
natural morphism
$$\wedge^pD^1\longrightarrow D^p$$
is an isomorphism and $(F_{\bullet}^{(p)}(D^1))$ is a minimal free
resolution of $D^p$, for every $p$, $p\leq n-1$.
 We have
$${\rm pdim}\,D^p=p,\,{\rm for}\,p\leq n-1,$$
$${\rm pdim}\,D^{n+1}=0,$$
and if ${\rm char}\,k$ does not divide $|{\cal A}|$, then
$${\rm pdim}\,D^n=n-1.$$
\end{thm}

\bibliographystyle{amsalpha}

\begin{thebibliography}{10}

\bibitem{bh} W. ~Bruns, J. ~Herzog,
         {\em Cohen-Macaulay rings},
         Cambridge studies in advanced mathematics 39, 1998.

\bibitem{d} P. ~Deligne,
        {\em Theorie de Hodge II},
        Inst. Hautes Etudes Sci. Publ. Math. \textbf{40} (1971), 5-58.

\bibitem{dk} I. ~Dolgachev, M. ~Kapranov,
        {\em Arrangements of hyperplanes and vector bundles on $P^n$},
        Duke Mathematical Journal \textbf{71} (1993), 633-664.

\bibitem{er} T. ~Edelman, V. ~Reiner,
        {\em A counterexample to Orlik's conjecture},
        Proceedings of the A.M.S. \textbf{118} (1993), 927-929.

\bibitem{e} D. ~Eisenbud,
        {\em Commutative Algebra with a view towards Algebraic Geometry},
        Graduate Texts in Mathematics, vol.~150,
        Springer-Verlag, Berlin-Heidelberg-New York, 1995.

\bibitem{esv} H.~Esnault, V. ~Schechtman, E. ~Viehweg,
        {\em Cohomology of local systems on the complement of hyperplanes},
        Inventiones Mathematicae \textbf{109} (1992), 557-561.

\bibitem{fu} W. ~Fulton,
        {\em Intersection Theory},
        Springer-Verlag, New York, 1984.

\bibitem{h} R. ~Hartshorne,
        {\em Algebraic Geometry}, Graduate Texts in Mathematics, vol. ~52,
        Springer-Verlag, Berlin-Heidelberg-New York, 1977.

\bibitem{l1} K. ~Lebelt,
        {\em Freie aufl\"osungen \"ausser potenzen},
        Maunscripta Mathematica \textbf{21} (1977), 341-355.

\bibitem{oss} C. ~Okonek, M. ~Schneider, H. ~Spindler,
        {\em Vector Bundles on Complex Projective Spaces}, Progress in
Mathematics,vol. ~3,
        Birkhauser, Boston, 1980.

\bibitem{ot} P. ~Orlik, H. ~Terao,
        {\em Arrangements of Hyperplanes}, Grundlehren Math. Wiss., Bd. ~300,
        Springer-Verlag, Berlin-Heidelberg-New York, 1992.

\bibitem{rt} L. ~Rose, H. ~Terao,
        {\em A free resolution for the module of logarithmic forms of a generic
arrangement},
        Journal of Algebra \textbf{136} (1991), 376--400.

\bibitem{sa} K. ~Saito,
        {\em Theory of logarithmic differential forms and logarithmic
vector fields},
        J. Fac. Sci. Univ. Toyko Sect. IA Math \textbf{27} (1981), 265-291.

\bibitem{s} H. ~Schenck,
        {\em A rank two vector bundle associated to a three arrangement,
and its Chern polynomial},
         Advances in Mathematics, to appear.

\bibitem{si} R. ~Silvotti,
        {\em On the Poincar\'e polynomial of a complement of hyperplanes},
        Mathematical Research Letters \textbf{4} (1997), 645--661.

\bibitem{st} L. ~Solomon, H. ~Terao,
        {\em A formula for the characteristic polynomial of an arrangement},
         Advances in Mathematics \textbf{64} (1987), 305-325.

\bibitem{t} H. ~Terao,
        {\em Generalized exponents of a free arrangement of hyperplanes and
Shepard-Todd-Brieskorn formula},
        Inventiones Mathematicae \textbf{63} (1981), 159--179.

\bibitem{y} S. ~Yuzvinsky,
        {\em A free resolution of the module of derivations for generic
arrangements},
        Journal of Algebra \textbf{136} (1991), 432-438.

\bibitem{y1} S. ~Yuzvinsky,
        {\em The first two obstructions to the freeness of arrangements},
        Transactions of the A.M.S. \textbf{335} (1993), 231-244.

\bibitem{y2} S. ~Yuzvinsky,
        {\em Free and locally free arrangements with a given intersection
lattice},
        Proceedings of the A.M.S. \textbf{118} (1993), 745-752.

\bibitem{z} G. ~Ziegler,
        {\em Combinatorial construction of logarithmic differential forms},
         Advances in Mathematics \textbf{76} (1989), 116-154.

\end{thebibliography}

\end{document}